\documentclass[a4paper,11pt]{article}
\usepackage[utf8]{inputenc}
\usepackage{textcomp}
\usepackage{graphicx}
\usepackage{graphics, amsmath,amsthm,amssymb}
\usepackage{bbm}
\usepackage{MnSymbol}
\DeclareMathAlphabet\mathbb{U}{msb}{m}{n}
\usepackage[T1]{fontenc}
\usepackage{float}
\usepackage{dsfont}
\newcommand{\blank}[1]{}
\newcommand{\vep}{\varepsilon}
\renewcommand{\P}{\mathbb{P}}
\newcommand{\N}{\mathbb{N}}
\newcommand{\Z}{\mathbb{Z}}

\newcommand{\E}{\mathbb{E}}
\newcommand{\ra}{\rightarrow}
\newcommand{\mc}[1]{\mathcal{#1}}
\newcommand\1{\mathrm{1}}
\newcommand{\indic}[1]{\1_{\{#1\}}}

\newtheorem{Lemma}{Lemma}[section]
\newtheorem{Proposition}[Lemma]{Proposition}

\newtheorem{Theorem}[Lemma]{Theorem}

\newtheorem{Definition}[Lemma]{Definition}

\newcommand{\Qed}{\qed \medskip}

\newcommand{\ceq}{\preccurlyeq}

\title{On strict monotonicity of the speed for\\ excited random walks in one dimension}
\author{Mark Holmes\footnote{Statistics Dept., The University of Auckland.  {\em m.holmes@auckland.ac.nz}}}
\begin{document}

\maketitle

\begin{abstract}
We give a ``direct'' coupling proof of strict monotonicity of the speed for 1-dimensional multi-excited random walks with positive speed.  This reproves (and extends) a recent result of Peterson without using branching processes.
\end{abstract}
\nocite{}
\section{Introduction and main results}
\label{sec:intro}
Multi-excited random walk is a nearest-neighbour self-interacting random walk model where, on the $k$-th departure from a site $x\in \Z^d$, the step distribution of the walker is given by $\beta(x,k,\bullet)$, supported on $\{\pm e_1,\pm e_2,\dots, \pm e_d\}$ (where $e_1,\dots, e_d$ denotes the canonical basis for $\Z^d$).  When $\{\{\beta(x,k,\bullet)\}_{k \in \N}\}_{x\in \Z^d}$ are i.i.d.~over $x$ the model is said to be a multi-excited random walk in i.i.d.~random cookie environment.  When $\beta(x,k,\bullet)=\beta(k,\bullet)$ for every $x$ (i.e.~the transition probabilities do not depend on the site of departure) the environment is said to be non-random.  The original model, introduced by Benjamini and Wilson \cite{BW03}, where $\beta(x,1,e)=(1+\beta e\cdot e_1)/(2d)$ and $\beta(x,k,e)=1/(2d)$ for all $e$ when $k>1$, has been extensively studied, and in particular it is known that the limiting velocity of the walker $\lim_{n \ra \infty}n^{-1}X_n=v(\beta)$ is (deterministic and) strictly positive when $d>1$ \cite{BR07,Kozm05,Kozm03} and strictly monotone in $\beta$ when $d>8$ \cite{HHmono} (see also \cite{PhamERW,Hol12,HSun12}).  It is believed that this monotonicity holds in dimensions $d>1$, but the absence of any natural coupling argument has been an obstacle to resolving this conjecture.  Note that the main idea behind proving monotonicity in high dimensions is to differentiate a formula for the speed provided by either an expansion \cite{HH12} or a Girsanov transformation, and it is currently not known how to adapt these arguments to low dimensions.  

In this paper we henceforth restrict ourselves to one dimension, where the law of the model is completely described by the quantities $\beta(x,k)\equiv \beta(x,k,e_1)$.  Here, Benjamini and Wilson's single cookie model is well understood (see e.g.~\cite{Dav99} and \cite{BW03}), but in the general setting of (multi-excited) random walks in random cookie environments (introduced by Zerner \cite{Zern05}) there are a number of interesting open problems (see the excellent recent survey of Kosygina and Zerner \cite{KZ13}).  

Amir, Berger, and Orenstein \cite{ABO13} (see also \cite[Theorem 4.1]{KZ13}) have recently confirmed that a deterministic speed exists when the environment is stationary and ergodic (in the $\Z$-shift) and elliptic ($\beta(x,k)\in (0,1)$). 
In i.i.d.~random environments, Basdevant and Singh, and Kosygina and Zerner, \cite{Zern05,BS08,BS08b,KZ08} have shown that under certain assumptions  the speed (exists, is deterministic and) is positive if and only if $\delta>2$, where
\begin{align}
\delta\equiv \E\left[\sum_{k=1}^\infty(2\beta(o,k)-1)\right].\label{delta}
\end{align}
Moreover the random walker is transient to the right if and only if $\delta>1$.  See \cite{KZ13} for a comprehensive survey of known results (up to 2012).

One such set of conditions under which the criterion \eqref{delta} for positive speeds and directional transience holds is when there exists $M\in \N$ such that $\beta(x,k)=p_k\in (0,1)$ for each $x$ (i.e.~non-random environment) and that $p_k=1/2$ for all $k\ge M$.  In this case the sum in \eqref{delta} is deterministic and finite.  This is the setting in which Peterson \cite{Pet14} proves a strict monotonicity property for the speed.  Before properly stating Peterson's results, we first recall the partial ordering on arrow systems introduced in \cite{HS_JCTA}.

An arrow system is an element $\mc{E}=\{\mc{E}(x,k)\}_{x \in \Z,k\in \N}$ of $\{-1, +1\}^{\Z\times \N}$.  An arrow system $\mc{E}$ defines (uniquely) a nearest neighbour walk $E=\{E_n\}_{n\in \Z_+}$ by setting $E_0=0$ and if $E_n=x$ and $\#\{m\le n:E_m=x\}=k$ then $E_{n+1}=E_n+\mc{E}(x,k)$.  For two arrow systems $\mc{L},\mc{R}$, we write $\mc{L}\ceq\mc{R}$ if $\sum_{j=1}^k\mc{L}(x,j)\le \sum_{j=1}^k\mc{R}(x,j)$ for each $x\in \Z$, $k\in \N$.  The following results are among those proved in \cite[Corollary 3.7, (3.4),Theorem 1.3]{HS_JCTA}.    
\begin{Theorem}
\label{thm:jcta}
If $\mc{L}\ceq \mc{R}$ then the corresponding walks $L,R$ satisfy:
\begin{itemize}
\item[(a)] For every $x>0$, $\inf\{n:R_n=x\}\le \inf\{n:L_n=x\}$.
\item[(b)] For every $n\in \Z_+$, $\max_{m\le n}L_m\le \max_{m\le n}R_m$ and $\min_{m\le n}L_m\le \min_{m\le n}R_m$.
\item[(c)] $\liminf_{n\ra \infty}L_n\le  \liminf_{n\ra \infty}R_n$.
\item[(d)] Let $a_n\le n$ be an increasing sequence, with $a_n\ra \infty$. If there exists $x\in \Z$ such that $R_n\ge x$ infinitely often
then $\limsup_{n \ra \infty}L_n/a_n\le \limsup_{n \ra \infty}R_n/a_n$.
\item[(e)] If $L_n \ra +\infty$ then $\#\{m:L_m=x\}\ge \#\{m:R_m=x\}$ for every $x$.
\end{itemize}
\end{Theorem}
In particular, (c) implies that if $L$ is transient to the right ($L_n \ra \infty$) then so is $R$, and $(d)$ implies that if the speeds $v(L)$ and $v(R)$  exist then $v(L)\le v(R)$.  These facts were used by the authors to prove monotonicity for the speed of excited random walks in one dimension, under some technical assumptions on the environments.  Note that examples are given in \cite{HS_JCTA} where $\mc{L}\ceq \mc{R}$ but $L_n>R_n$ for some $n$ etc.   

Peterson \cite{Pet14} strengthened these results to strict monotonicity of the speed under much stronger assumptions.  To be precise, Peterson considers the situation where there are two elliptic non-random i.i.d.~cookie environments $\vec{p}\in (0,1)^{\N}$ and $\vec{q}\in (0,1)^{\N}$, such that $\vec{p}\prec\vec{q}$ where: 
\begin{Definition}
\label{def:peterson}
$\vec{p}\prec \vec{q}$ if there exists a coupling $\P$ of $(\vec{Y},\vec{Z})$ with $\vec{Y} = \{Y_i\}_{i \in \N}$
and $\vec{Z} = \prec\{Z_i\}_{i \in\N}$ such that $\{Y_i\}_{i \in \N}$ are independent random variables with $Y_j\sim $Ber$(p_j)$,
$\{Z_i\}_{i \in\N}$ are independent with $Z_j\sim$ Ber$(q_j)$, and moreover
\begin{align}
\P\Big(\sum_{j=1}^mY_j\le \sum_{j=1}^m Z_j\Big)=1 \quad \textrm{ for every }m\in \N,\label{weakassumption}\\
\P\Big(\sum_{j=1}^mY_j<\sum_{j=1}^m Z_j\Big)>0 \quad \textrm{ for some }m\in \N.\label{strongassumption}
\end{align}
\end{Definition}
Swapping $p_i$ and $p_j$ (i.e.~setting $p'_j=p_i$ and $p'_i=p_j$) is said to be {\em favorable} if $i<j$ and $p_j>p_i$.  
Given an environment $\vec{p}$, one should think of generating environments $\vec{q}\succ \vec{p}$ in two generic ways:  increasing $p_i$ for one or more $i$, or making favorable swaps of $p_i$.

Peterson then proves the following by appealing to a connection with branching processes with migration (which has become a standard tool in proving ballisticity and related properties for multi-excited random walks on integers \cite{KZ13}), and comments that ``proving these results via a direct coupling of excited random walks fails''.
\begin{Theorem}[Theorems 1.7 and 1.8 of \cite{Pet14}]
\label{thm:peterson}
Let $\vec{p},\vec{q}\in (0,1)^{\N}$.  Suppose that 
\begin{itemize}
\item[(i)] $\vec{p}\prec\vec{q}$, and 
\item[(ii)] there exists $M\in \N$ such that $p_k=q_k=1/2$ for all $k>M$. 
\end{itemize}
Then either $v(\vec{p}) = v(\vec{q}) = 0$ or $ v(\vec{p}) < v(\vec{q})$.
Moreover, if the $\vec{p}$ walk is transient to the right (a.s.) then $\P_{\vec{p}}(X_n > 0, \forall n > 0) < \P_{\vec{q}}(X_n > 0, \forall n > 0)$.  
\end{Theorem}
In this paper we reprove this result without assuming (ii).  In other words, we prove:
\begin{Theorem}
\label{thm:main}
Let $\vec{p},\vec{q}\in (0,1)^{\N}$ be such that $\vec{p}\prec\vec{q}$.
Then the conclusions of Theorem \ref{thm:peterson} hold.
\end{Theorem}
The existence of (deterministic) $v(\vec{p})$ and $v(\vec{q})$ is known to hold in this elliptic (and deterministic) setting by \cite{ABO13}.  Our proof is via a coupling of arrow systems, which might be considered a ``direct'' coupling of excited random walks, in the sense that it does not use any facts about branching processes.  We note that our methods and results can be extended to more general cookie environments (see \cite[Section 5]{HS_JCTA}, and also Section \ref{sec:general}) at the cost of much more cumbersome notation.  

One class of environments to which condition (ii) does not hold yet our results apply, is the so-called {\em periodic} environments, as studied by Kozma, Orenshtein, and Shinkar in \cite{periodic}.  In this setting the finite sequence of cookies $p_1,\dots,p_M$ is repeated indefinitely, i.e.~$p_{i+M}=p_i$ for every $i\in \N$. Let $\bar{p}=M^{-1}\sum_{i=1}^M p_i$ and 
\begin{align*}
\theta(p_1,\dots,p_M)=\frac{\sum_{i=1}^M (1-p_i)\sum_{j=1}^i(2p_j-1)}{4\sum_{l=1}^Mp_l(1-p_l)}.
\end{align*}
It is shown in \cite{periodic} that the walk is transient to the right if either $\bar{p}>1/2$ or $\bar{p}=1/2$ and $\theta(p_1,\dots,p_M)>1$.  Although no criterion for positivity of the speed is given in \cite{periodic}, it is certainly the case that the speed is positive when $\bar{p}>1/2$, and also that the speed will depend on the order of the cookies.  We expect that one can also obtain positive speeds when $\bar{p}=1/2$, depending e.g.~on the order of the cookies.

Theorem \ref{thm:main} applied in the context of periodic cookie environments says that if $\vec{p}$ is periodic and induces a positive speed then any environment $\vec{q}$ created by favourable swaps has greater speed.  This of course includes the cases where $\vec{q}$ is also periodic.


\section{Proof of Theorem \ref{thm:main}}
\label{sec:proof}
Our first result is an elementary one connecting arrow systems and excited random walks.
\begin{Lemma}
\label{lem:space}
If $\vec{p}\prec \vec{q}$ then there exists a probability space $(\Omega, \mc{F},\P)$ on which there are (random) arrow systems $\mc{L}\ceq\mc{R}$ whose corresponding walks $L,R$ have the laws of $\vec{p}$ and $\vec{q}$ excited random walks respectively.
\end{Lemma}
\proof By assumption there exists a probability space $(\Omega, \mc{F},\P)$ on which we have random sequences $(\vec{Y},\vec{Z})$ satisfying the properties in Definition \ref{def:peterson}.  By considering product spaces one can easily embed this into a larger probability space  where we have $\{(\vec{Y},\vec{Z})_x\}_{x\in \Z^d}$ each with the same law and independent over $x$.  
Now define $\mc{L}(x,k)=2Y_{x,k}-1$ and $\mc{R}(x,k)=2Z_{x,k}-1$.  It is easy to see that $\mc{L}, \mc{R}$ and $L,R$ have the desired properties.\Qed

Standard proofs of the existence of a deterministic speed are based on the notion of regeneration points and times.  A point $x\in \Z_+$ is a regeneration point for a (right-transient) walk $X$ if $X_n\ge x$ for all $n\ge \inf\{m\in \N:X_m=x\}$, in other words if the walker never returns to the left of $x$ after it reaches $x$.  Let $\mc{D}(X)$ denote the set of  (non-negative) regeneration points of $X$.  If $\mc{E}$ is an arrow system, let $\mc{D}(\mc{E})=\mc{D}(E)$ denote the regeneration points of the walk $E$ defined from it.  

The following is now a standard result in the literature, see e.g.~\cite[Lemma 3.17 (also Lemma 4.5, Theorem 4.6)]{KZ13}.
\begin{Lemma}
\label{lem:standard}
For an i.i.d.~elliptic environment, if $\P(X_n\ra +\infty)>0$ then for every $x\in \Z_+$, $\P(x\in \mc{D}(\vec{p}))=\vep_{\vec{p}}>0$.
\end{Lemma}

The following elementary lemma is one of the main facts that lets us upgrade monotonicity as proved in \cite{HS_JCTA} to strict monotonicity (under stronger assumptions such as those of Theorem \ref{thm:main}).  In this lemma, $T_E(x)$ denotes the first hitting time of level $x>0$ by the walk $E$, and $\mc{D}(L)=\{D_1(L),D_2(L),\dots\}$ where $0\le D_i(L)<D_{i+1}(L)$ for each $i$.
\begin{Lemma}
\label{lem:regen}
If $\mc{L}\ceq\mc{R}$ then $\mc{D}(\mc{L})\subset \mc{D}(\mc{R})$ on the probability space in Lemma \ref{lem:space}.  Moreover $T_R(D_{k+1}(L))-T_R(D_{k}(L))\le T_L(D_{k+1}(L))-T_L(D_{k}(L))$ for every $k\in \N$.
\end{Lemma}
\proof Theorem \ref{thm:jcta}(b) proves that if $0\in \mc{D}(\mc{L})$ then $0\in \mc{D}(\mc{R})$.  If $L$ reaches level $x$ then Theorem \ref{thm:jcta}(b) proves that $R$ does too.  Both claims then follow immediately since the property $\mc{L}\ceq\mc{R}$ is translation invariant, and $x\in \mc{D}$ does not depend on the arrow system to the left of $x$.\Qed

Let $T(E)=\inf\{n\in \N:X_n\in \mc{D}(E)\setminus 0\}$ denote the first hitting time of the first regeneration level to the right of the starting point, and for $\vec{p}\in (0,1)^{\N}$, let $T(\vec{p})$ denote the corresponding hitting time for an excited random walk with cookie environment $\vec{p}$.  The following result connects the speeds of (a.s.~right-transient) excited random walks $L\ceq R$ given by \cite[Lemma 4.5, Theorem 4.6]{KZ13} according to their {\em mutual} regeneration levels.  Note that coupling of regenerations is not new (see e.g.\cite{BAFGH12}, where the notion of super-regeneration times are introduced), but given the importance in the present context, we include a proof of the following result.

\begin{Proposition}
\label{prp:regenerationtimes}
Assume that the excited random walk in environment $\vec{p}\in (0,1)^{\N}$ is almost surely transient to the right.  If $\vec{p}\prec \vec{q}$ then there exists a probability space on which
\begin{align}
v(\vec{p})=&\frac{\E[X_{T(\vec{p})}\indic{0 \in \mc{D}(\vec{p})}]}{\E[T(\vec{p})\indic{0 \in \mc{D}(\vec{p})}]}, \quad \text{and}\label{speedp}\\
 v(\vec{q})=&\frac{\E[X_{T(\vec{p})}\indic{0 \in \mc{D}(\vec{p})}]}{\E[T_{\vec{q}}(\vec{p})\indic{0 \in \mc{D}(\vec{p})}]},\label{speedq}
\end{align}
where $T_{\vec{q}}(\vec{p})$ denotes the hitting time of $X_{T(\vec{p})}$ by the $\vec{q}$-walk.
\end{Proposition}
\proof Construct the probability space given by Lemma \ref{lem:space}, on which $L$ and $R$ are $\vec{p}$ and $\vec{q}$ excited random walks defined from $\mc{L}$ and $\mc{R}$ respectively.

Since $L$ is almost surely transient to the right (by assumption), 
the first claim is then a simple consequence of \cite[Lemma 3.17, Lemma 4.5, Theorem 4.6]{KZ13}, using the fact that the chunks of arrow environment seen by the $L$-walker in between regeneration times (times between hitting of regeneration levels) are i.i.d.  By Theorem \ref{thm:jcta}(c), $R$ is also transient to the right (a.s.), and by Lemma \ref{lem:regen} we have $\mc{D}(\mc{L})\subset \mc{D}(\mc{R})$ on this probability space.  Since the arrow environments seen by the walk $R$ in between hitting regeneration levels of $L$ are i.i.d., the second claim follows in exactly the same way. 
\Qed

\noindent {\em Proof of Theorem \ref{thm:main}.}  We continue to work in the probability space of Lemma \ref{lem:space}.  Without loss of generality we may assume that $v(\vec{p})\ge 0$.  

On the event $\{L_n \nrightarrow +\infty\}$, we have $L_n\le x$ infinitely often for some $x$, whence the $\liminf$ analogue of Theorem \ref{thm:jcta}(d) applies and tells us that $0=\liminf n^{-1}L_n\le \liminf n^{-1}R_n$.  

On the event $\{L_n \ra +\infty\}$, regeneration levels exist (a.s.) and Lemma \ref{lem:regen} tells us that $\mc{D}(\mc{L})\subset \mc{D}(\mc{R})$ and the inter-hitting times of the regeneration levels satisfy $T_R(D_{k+1}(L))-T_{R}(D_{k}(L))\le T_L(D_{k+1}(L))-T_L(D_{k}(L))$ for every $D_{k}(L)<D_{k+1}(L)\in \mc{D}(L)$.  It follows that on the event $L_n \ra +\infty$: the speeds $v_L=\lim_{n \ra \infty} n^{-1}L_n$ and $v_R=\lim_{n \ra \infty} n^{-1}R_n$ (exist and) satisfy $v_L\le v_R$ (a.s.).  This gives (non-strict) monotonicity of the speeds on the event $L_n \ra +\infty$.

Combining the arguments for the two events above proves the claimed monotonicity of the speed when $v(\vec{p})=0$.  

Suppose now that $v(\vec{p})>0$.  Then $L_n \ra +\infty$ almost surely so $v(\vec{q})\ge v(\vec{p})$ by Theorem \ref{thm:jcta}(d), and Proposition \ref{prp:regenerationtimes} gives us formulas for the speeds.
In particular, since $v(\vec{p})>0$ we have $\E[T(\vec{p})\indic{0 \in \mc{D}(\vec{p})}]<\infty$.  To verify strict monotonicity of the speed it therefore remains to show that 
\begin{align}
\P(T(\vec{p})\indic{0 \in \mc{D}(\vec{p})}>T_{\vec{q}}(\vec{p})\indic{0 \in \mc{D}(\vec{p})})>0.\label{strict}
\end{align}
By Lemma \ref{lem:standard}, $\vep_{\vec{p}}\equiv \P(2 \in \mc{D}(L))>0$.  Note that this is independent of the environment before level $2$ (assuming one reaches level 2).  Let $m_0=\inf\{m\in \N:\P(\sum_{j=1}^mY_j<\sum_{j=1}^m Z_j)>0\}$, which is finite since $\vec{p}\prec \vec{q}$.  Since $p_k\in (0,1)$ for each $k$, with positive probability on our probability space, all of the following hold simultaneously
\begin{enumerate}
\item[(i)] $Y_{j,0}=Z_{j,0}=1$ for all $j\le m_0+1$,
\item[(ii)] $Y_{j,1}=0$ for all $j\le m_0$, $Z_{m_0,1}=1$,
\item[(iii)] $2\in \mc{D}(L)\subset \mc{D}(R)$.
\end{enumerate}
It is then easy to see that on this event $0 \in \mc{D}(\vec{p})\subset \mc{D}(\vec{q})$, and that $T(\vec{p})=\inf\{n:L_n=2\}=2m_0+2>2m_0\ge \inf\{n:R_n=2\}=T_{\vec{q}}(\vec{p})$, as required.

The last claim of the theorem is that $\P(0 \in \mc{D}(\vec{q}))>\P(0 \in \mc{D}(\vec{p}))$ when the latter is positive.  This follows easily by a local construction similar to that above.
\Qed

\section{Further generalisations}
\label{sec:general}
Note that all we have really needed here to make the proof work is: 
\begin{itemize}
\item[(1)] there is a coupling of the walks in terms of arrow systems $\mc{L}\ceq \mc{R}$, 
\item[(2)] the environment is i.i.d., and regeneration points exist when we have transience to the right (so that speeds exist as well, otherwise we should not compare them), and 
\item[(3)] under this coupling, with positive probability the hitting time of the first regeneration level is smaller for $R$ than $L$ (as in \eqref{strict}).  Assuming that with positive probability the two environments actually differ, this is easy to achieve under the assumption of ellipticity.
\end{itemize}
In fact the i.i.d.~assumption is really only needed to give comparable formulas  for the speeds as in \eqref{speedp} and \eqref{speedq}.  So if for example two walks can be coupled via  arrow systems arising from stationary ergodic environments such that formulas for the speeds can be given in terms of regeneration levels, then the same proof works in this case.  Similarly the method used here remains valid if one has for example a pair of right transient environments alternating (or more generally cyclic) in $x \in \Z$  as
$\vec{p}_{\textrm{even}},\vec{p}_{\textrm{odd}}$ and $\vec{q}_{\textrm{even}},\vec{q}_{\textrm{odd}}$, with $\vec{p}_{\textrm{even}}\prec\vec{q}_{\textrm{even}}$ and $\vec{p}_{\textrm{odd}}\prec\vec{q}_{\textrm{odd}}$ (only one would need to be strict).  In order for the strict monotonicity proof to go through one would need formulas for the speeds, so in particular an external input of existence of speeds and e.g.~``even regeneration levels''.

\section*{Acknowledgements}
This work is supported in part by the Marsden Fund, administered by RSNZ.  The author thanks Tom Salisbury and Gady Kozma for helpful discussions.  

\bibliographystyle{plain}

\end{document}